\documentclass[aps,pra,11pt]{revtex4-1}
\usepackage{amsmath,amsfonts,amssymb}
\usepackage{latexsym}

\usepackage[active]{srcltx} % para asociar búsqueda desde el viewer de dvi al editor (Okular a Kile) con Shift-Click

\newtheorem{theorem}{Theorem}[section]
\newtheorem{definition}[theorem]{Definition}
\newtheorem{example}[theorem]{Example}

\newtheorem{pro}[theorem]{Proposition}

\def\G{\Gamma}

\def\RR{\mathbb{R}}
\def\NN{\mathbb{N}}

\def\CC{\mathbb{C}}

\def\ket{\rangle}
\def\bra{\langle}
\def\d2z{d^{2}z}
\def\beq{\begin{equation}}
\def\eeq{\end{equation}}
\def\beqn{\begin{eqnarray}}
\def\eeqn{\end{eqnarray}}

\begin{document}

\title{About radial Toeplitz operators on Segal-Bargmann and $l^2$ spaces.}

\author{Romina A.\ Ram\'{\i}rez}
\affiliation{
Departamento de Matem{\'a}tica, Facultad de Ciencias Exactas,
Universidad Nacional de La Plata, Argentina.
}
\email{romina@mate.unlp.edu.ar}

\author{Gerardo L.\ Rossini}
\affiliation{
IFLP, CONICET - Departamento de F\'{\i}sica and  Departamento de Matem{\'a}tica, Facultad de Ciencias Exactas,
Universidad Nacional de La Plata, Argentina.
}
\email{rossini@fisica.unlp.edu.ar}

\author{Marcela Sanmartino}
\affiliation{Departamento de Matem{\'a}tica, Facultad de Ciencias Exactas,
Universidad Nacional de La Plata, Argentina.}
\email{tatu@mate.unlp.edu.ar}
%\date{\today}

%\subjclass{Primary 47B35; Secondary 81Q99}

%\keywords{Toeplitz operators; radial symbols}

%\date{}

\begin{abstract}

We discuss Toeplitz operators on the Segal-Bargmann space as functional realizations of 
anti-Wick operators on the Fock space.
In the special case of radial symbols we exploit the isometric mapping
between the Segal-Bargmann space and $l^2$ complex sequences in order to
establish conditions such that an equivalence between Toeplitz operators
and diagonal operators on $l^2$ holds.
We also analyze the inverse problem of mapping diagonal operators on $l^2$
into Toeplitz form.
The composition problem of Toeplitz operators with radial symbols is reviewed as an application. 
Our notation and  basic examples make contact with Quantum Mechanics literature.

\end{abstract}

\maketitle
%\begin{keyword}
%{Toeplitz operators}
%
%\end{keyword}
%\noindent{\bf 2000 AMS Subject Classification}: Primary 47B35;
%Secondary 47G30

%\newpage
\section{\textbf{Introduction}}
\label{intro}

Toeplitz operators were introduced in physics by Berezin \cite{B1} \cite{B2} \cite{Berezin} in the
context of quantization procedures, {\em i.e.} in the association of a classical function in phase
space and a quantum observable.
In this sense, the classical function is called a symbol for the operator.
Indeed, there are different ways to quantize, all of them consistent
with the probabilistic interpretation of quantum mechanics.
The key difference between them arises from non-commutativity of quantum observables, in contrast
with commutativity of their classical counterparts.
To be more specific, in a Fock space there exist annihilation and creation operators
$\widehat{a}$, $\widehat{a}^{*}$ which do not commute;
the symbolic calculus associated to the product order $\widehat{a}^{*}\widehat{a}$ leads to Wick
symbols, while that
associated to the product order $\widehat{a}\widehat{a}^{*}$ leads to anti-Wick symbols.
In this context, Toeplitz operators arise as the Segal-Bargmann space realizations of 
Fock space operators having an associated anti-Wick symbol (formal definitions are given below).

More generally, given a symbol $\varphi(z,\overline{z})$, the corresponding Toeplitz operator 
is defined  as the Bargmann projection of the pointwise multiplication 
of a function in the Segal-Bargmann space with the symbol.
In this sense, the operator has a natural domain not covering the Segal-Bargmann space, 
but restricted to functions such that the projection is well defined. This restriction poses a problem 
on various properties of Toeplitz operators, such as composition.

Toeplitz operators are also defined in the Bergman space of analytical functions on the unit disk,
where most of the theoretical results have been developed (see for instance \cite{Bergman}). 
They have also been object of study in different 
disciplines. 
In partial differential equations, these operators and their
adjoints play an important role in extending known results in the
space of entire functions to the context of Segal-Bargmann
spaces (see for instance \cite{Ci}, \cite{Scha}), \cite{Ja1}, \cite{Janas1994}).
They have also been extensively studied as an efficient
mathematical tool in signal analysis (\cite{C-G}, \cite{C-R},
\cite{G}).

It is known that the composition of Toeplitz operators is in
general not closed, in the sense that the it may be well defined but the resulting operator
is not a Toeplitz one.
The problem of how to define the class of the symbols where the
composition of the corresponding Toeplitz operators is closed is still open.
Some authors have addressed this issue:
\cite{Coburn} presents some classes of operators where the problem is solved, while
recent extensions are presented in \cite{Bauer}.
In contrast, in \cite{A-M}, \cite{C-G}, \cite{C-R} and \cite{L} the composition of certain Toeplitz
operators has been shown to be expressed as a Toeplitz operator plus a remainder term.

This article is intended to provide a unified discussion of operators with radial anti-Wick symbols. 
Following insights from Quantum Mechanics, we start from the abstract Fock space formulation and 
discuss concrete functional realizations both in Segal-Bargmann  and $l^2$ spaces.
We then apply the results to the related composition problem.

The present work is organized as follows.
Section \ref{II} and the Appendix contain preliminary material about Fock space and its functional realizations, in 
particular the isometry between the Segal-Bargmann and $l^2$ spaces.
In Section \ref{III} precise definitions of anti-Wick and Toeplitz operators are given, stressing the class of 
symbols and operator domains considered by different authors.
Section \ref{IV} addresses the case of Toeplitz operators with radial symbols.
Under the isometry between the Segal-Bargmann space and $l^2$ sequences, some of them are unitarily
equivalent to diagonal operators on $l^2$.
While this drastic simplification allows for a very simple analysis of their properties \cite{G-V}, one should 
notice that not any Toeplitz operator with radial symbol can be treated in this way. We provide 
sufficient and necessary conditions for this equivalence to hold.
Section \ref{V} aims to investigate some natural questions related to the the inverse problem of existence and uniqueness 
of an anti-Wick symbol and Toeplitz operator for a given diagonal operator on $l^2$. 
While the general problem is hard analyze, 
we present a family of such operators for which the symbols 
can be explicitly constructed and study, on this set, sufficient and necessary conditions for 
unitarily equivalence between diagonal operators on $l^2$ and Toeplitz operators. 
The studied family is large enough to illustrate accomplishment or not of this conditions.
Section \ref{VI} discusses the composition problem, addressing to the large gap between
classes of Toeplitz operators for which the composition is known to be closed and counter examples where
the composition is a well defined but not a Toeplitz operator \cite{Coburn}.
We resort to our results in Section \ref{V} to present novel results on composition, and relate them to known 
positive results and counter-examples.
%

%%%%%%%%%%%%%%%%%%%%%%%%%%%%%%%%%%%%%%%%%%%%%%%%%%%%%%%%%%%%%%%%%%%%%%%%%%%%%%%%%%%%%%%%%%%%%%%%%%%%%%%%%%%%%%%%%%%%%%%%%%

\section{Functional realizations of the Fock space}\label{II}

An abstract Fock space ${\mathcal F}$ is a Hilbert space of vectors, denoted by $|\psi\ket$  in Dirac's notation,
in which there exist an
operator $\widehat{a}$ and its adjoint $\widehat{a}^{*}$, called annihilation
and creation operators, satisfying the canonical commutation
rules  $[\widehat{a},\widehat{a}^{*}]=I$, where $I$ is the
identity operator on ${\mathcal F}$.
There also exists a vector $|0\ket$
(called \emph{ vacuum vector}) annihilated by $\widehat{a}$, such that the system
\beq
\left\{
|n\ket = \frac{(\widehat{a}^{*})^n|0\ket }{\sqrt{n!}}
\right\}_{n \in \NN}
\label{Fockbasis}
\eeq
is complete and orthonormal in ${\mathcal F}$ (the canonical Fock space basis, 
known in physics as the {\em occupation number} basis).
Then ${\mathcal F}$ is the space of linear combinations 
\beq
| \psi \ket = \sum_{n \in \NN} \psi_n | n \ket
\label{ket}
\eeq
with complex coefficients and finite norm. The inner product is noted as $\bra\phi|\psi\ket$, being antilinear in the
left vector and linear in the right one.

Functional realizations of the Fock space are obtained by projecting vectors $|\psi \ket$ on complete sets
labeled by either a discrete or continuous variable.
Classical examples of Fock space realizations are:

\vspace{5mm}
$\mathbf {l^2}$: projection on the orthonormal basis (\ref{Fockbasis}) 
gives simply  $ \bra n | \psi \ket = \psi_n $. The finite norm condition on $ | \psi \ket$ 
is equivalent to $\sum_{n} |\psi_n|^2 < \infty$.
Then the vector $| \psi \ket$ in ${\mathcal F}$ is realized by the sequence $\{\psi_n\}$ in $l^2$, 
the linear space of square summable complex sequences. The inner product in $l^2$ is realized by
$$
(\{\psi'_n\}, \{\psi_n\}) := \sum_{n\in\NN} \overline{\psi'_n}\psi_n =  \bra \psi' | \psi\ket   \, .
$$

\vspace{5mm}
$\mathbf {L^2(\RR, dx)}$: consider the set
$\left\{|x\ket\right\}$  containing the generalized eigenvectors of the (dimensionless) position operator
$\widehat{x}=(\widehat{a}^{*} + \widehat{a})/\sqrt{2}$, with continuum spectrum $\RR$ 
(see for instance \cite{Sakurai}). 
Projection on this set provides a space of complex valued functions $\psi(x)=\bra x | \psi \ket$.
The finite norm of $ | \psi \ket$ in ${\mathcal F}$ implies that $\psi(x)$ is
square integrable on $\RR$ with the Lebesgue measure $dx$. The inner product here is realized by
$$
(\psi',\psi) := \int_\RR \overline{\psi'(x)} \psi(x) dx = \bra \psi' | \psi\ket   \, .
$$
Then this realization leads to the linear space $L^2(\RR, dx)$. 
This is the most usual representation in physics,
known as coordinate representation, while its elements are called wave functions.
Notice that $|x\ket$ is not a vector in ${\mathcal F}$, but in the corresponding Rigged Hilbert space.

\vspace{5mm}
$\mathbf {F^{2}(\CC, d\mu)}$: Segal-Bargmann space: consider the set
\beq
\left\{|z\ket : \widehat{a}|z\ket = z |z\ket \right\}
\label{CS}
\eeq
containing the normalized eigenvectors of $\widehat{a}$, with continuum spectrum $\CC$ (see for instance \cite{Perelomov}).
The vectors $ |z\ket $ are called coherent states (or Poisson vectors) and form an overcomplete set in ${\mathcal F}$
(see the Appendix for technical details).
Projection on this set provides complex valued functions $\bra \overline{z} | \psi \ket$.
It is convenient to write
\beq
\psi(z) = \bra \overline{z} | \psi \ket  e^{\frac{1}{2}|z|^2}, 
\label{SBsymbol}
\eeq
called Segal-Bargmann symbol for $| \psi\ket \in {\mathcal F} $ in the following,
because the finite norm of $ | \psi \ket$ in ${\mathcal F}$ implies $\psi(z)$ is an entire function. 
Moreover, $\psi(z)$ is square integrable with the Gaussian measure  
$d\mu(z)=\frac{1}{\pi}e^{-|z|^{2}}\d2z$ 
($\d2z$ being the Lebesgue translationally invariant measure on $\CC$).
This realization, called the Segal-Bargmann space $F^{2}(\CC, d\mu)$, is then
the subset of {\em entire} functions in $L^{2}(\mathbb{C}, d\mu)$,
the Hilbert space of square integrable functions on $\CC$ with the Gaussian
measure $d\mu(z)$.
The inner product in $F^{2}(\CC, d\mu)$ is realized by
$$
(\psi',\psi) :=  \int_\CC \overline{\psi'(z)} \psi(z) d\mu(z) = \bra\psi'|\psi\ket  \, .
$$

\vspace{1cm}
The Segal-Bargmann functional realization $F^{2}(\CC, d\mu)$, where Toeplitz operators are defined, 
has as important advantages the connection with the powerful theory of analytic functions,
and the fact that vectors  $ |z\ket \in {\mathcal F}$ (in contrast with  $ |x\ket \notin {\mathcal F}$). 
However, the set of coherent states is not an orthogonal basis but an overcomplete set in ${\mathcal F}$.
Completeness means that vectors in ${\mathcal F}$ can be written as 
\beq
\int_\CC d\mu(z) f(z,\overline{z})\, |\overline{z}\ket.
\label{general-ket}
\eeq
Indeed, one can show that such a linear combination is a vector in ${\mathcal F}$
if and only if   $f(z,\overline{z}) \in L^2(\CC, d\mu)$.  
Overcompleteness means that the function $f(z,\overline{z})$ is not unique for a given vector $| \psi \ket\in {\mathcal F}$;
all the functions in $L^2(\CC, d\mu)$ expanding the same vector define equivalence classes, but there exists 
a unique entire function $\psi(z) \in  F^2(\CC, d\mu)$ representing each class, 
exhibiting the isomorphism between ${\mathcal F}$ and  $F^2(\CC, d\mu)$.
Given a vector $| \psi \ket$, such a function 
is just the Segal-Bargmann symbol $\psi(z)$ of the vector, defined in eq. (\ref{SBsymbol}).
These properties are easily derived using the well known Bargmann projection,
defined in $L^{2}(\CC, d\mu)$ as follows: given $g \in L^{2}(\CC, d\mu)$, let
\beq
K_{\overline{z}}(w) := e^{\overline{z} w}
\label{kernel}
\eeq
and
\beq
(P ~ g)(z) = (K_{\overline{z}}, g),
\label{projector}
\eeq
with the inner product in $L^{2}(\CC, d\mu)$. It holds that 
$$
P:L^{2}(\CC, d\mu)\rightarrow {F}^{2}(\CC, d\mu) ;
$$
this feature is essential to the definition of Toeplitz operators given below in eq. (\ref{Toeplitz}).

From the description above, there exist natural unitary isomorphisms between $l^2$, $L^2(\RR, dx)$ and $F^{2}(\CC, d\mu)$.
Indeed, since the Fock space is a unitary irreducible representation of the Heisenberg-Weyl group
\cite{Perelomov}, these isomorphisms are just examples of the celebrated theorem by Stone and von Neumann \cite{Stone,vonNeumann} stating that any two unitary irreducible
representations of the Heisenberg-Weyl group are unitarily equivalent.
It is then convenient to use the different realizations for the purposes to which they are most suited.
All the results proved for one realization, if formulated only in terms of a vacuum vector and
annihilation an creation operators, are then valid for the others.

Let us recall that the natural isomorphism between $F^2(\CC, d\mu)$ and $l^2$, 
that will be used throughout the present work, is derived from the relation 
\beq
\bra \overline{z}|n\ket = e^{-\frac{1}{2}|z|^2} \frac{z^{n}}{\sqrt{n!}}
\label{zn}
\eeq
(see the Appendix) between vectors in the canonical basis (\ref{Fockbasis}) and the coherent states set (\ref{CS}).
A vector $|\psi\ket$, as given in (\ref{ket}), is realized in $l^2$ by $\{\psi_n\}$ and in $F^2(\CC, d\mu)$ by
$ \psi(z) = \sum_{n\in \NN} \psi_n \frac{z^{n}}{\sqrt{n!}}$.

We summarize this relations in the following:
\begin{pro}
Let $U:F^2(\CC, d\mu) \to l^2$ be defined by 
\beq
U\psi = \{ \psi^{(n)}(0)/\sqrt{n!} \} .
\label{U}
\eeq
Then $U$ is a unitary isomorphism between $F^2(\CC, d\mu)$ and $l^2$.
\label{prop_U}
\end{pro}

{\em Proof}: Considering the Maclaurin expansion of $\psi \in F^2(\CC, d\mu)$, it is straightforward to compute
$$
(\psi,\psi) = \sum_{n\in\NN} \frac{|\psi^{(n)}(0)|^2}{n!},
$$
showing that $\{\psi_n\} \equiv \{ \psi^{(n)}(0)/\sqrt{n!} \} \in l^2$. 
Because of uniqueness of the expansion coefficients, $U^{-1}: l^2 \to F^2(\CC, d\mu)$ exists and is given by
$$
(U^{-1} \{\psi_n\})(z) = \sum_{n\in\NN}  \frac{\psi_n}{\sqrt{n!}} z^n,
$$
with the series converging in $\CC$.

Finally, it is easy to compute for any $\psi', \, \psi \in F^2(\CC, d\mu)$ that
$$
(U\psi',U\psi) =  \sum_{n\in\NN} \overline{\psi'_n} \psi_n = (\{\psi'_n\}, \{\psi_n\})
$$
showing that $U$ is unitary.

%%%%%%%%%%%%%%%%%%%%%%%%%%%%%%%%%%%%%%%%%%%%%%%%%%%%%%%%%%%%%%%%%%%%%%%%%%%%%%%%%%%%%%%%%%%%%%%%%%%%%%%%%%%%%%%%%%%%%%%%%%

\section{Anti-Wick and Toeplitz operators}\label{III}

We now introduce a class of integral operators on ${\mathcal F}$ formally given by a ``diagonal'' expression
in the coherent states set. Namely, given a measurable function $\varphi(w,\overline{w})$, not necessarily analytic, let
\beq
A_\varphi = \int_\CC  \frac{d^2w}{\pi} \, |\overline{w}\ket \varphi(w,\overline{w}) \bra \overline{w} |
\label{antiWick}
\eeq
where the notation $\bra \overline{w} |$  stands for the linear form $\bra\overline{w} | \cdotp \ket$ on ${\mathcal F}$ 
and the integral is understood in the weak sense. The function $\varphi(w,\overline{w})$ is known
as the anti-Wick or contravariant symbol  of $A_\varphi$ \cite{B1}, and we refer to $A_\varphi$
as an anti-Wick operator. A simple example is the (dimensionless) Hamiltonian operator for a harmonic oscillator,
$\widehat{H}= \widehat{a}^{*} \widehat{a} +1/2$, whose anti-Wick symbol results $H(w,\overline{w}) = w \overline{w} -1/2$; 
such an operator form can indeed be obtained for any polynomial operator in $\widehat{a}$ and  $\widehat{a}^{*}$  \cite{Berezin}.

For a given symbol, the formal expression (\ref{antiWick}) defines an operator $A_\varphi: D \subseteq {\mathcal F} \to {\mathcal F} $
with non trivial domain $D$ when the symbol is suitably restricted. For instance,
if, for some $r>2$,
$\int_\CC |\varphi(w,\overline{w})|^r \exp(-|w|^2)  \frac{d^2w}{\pi} < \infty$ , then
$A_\varphi$ is well defined on a dense domain in ${\mathcal F}$ \cite{Perelomov}.

Being written in terms of coherent state vectors, the most natural realization of $A_\varphi$ is that on ${F}^{2}(\CC, d\mu)$.
Projecting $A_\varphi | \psi \ket$ on $| \overline{z}\ket$, 
\beq
\bra \overline{z} |A_\varphi | \psi \ket
 =
\int_\CC  \frac{d^2w}{\pi} \bra \overline{z} |\overline{w}\ket \varphi(w,\overline{w})\bra \overline{w} |\psi \ket
\label{realizeF2}
\eeq
leads (see the Appendix) to the following expression for its Segal-Bargmann symbol: 
\beq
%(T_\varphi\psi)(z) = 
\int_\CC d\mu(w)\overline{K_{\overline{z}}(w)} \varphi(w,\overline{w}) \psi(w).
\label{Toeplitz}
\eeq
Let us call
$$
T_\varphi: D \subseteq F^2(\CC, d\mu) \to  F^2(\CC, d\mu)
$$ 
the realization of $A_\varphi$ acting on ${F}^{2}(\CC, d\mu)$.
As the integral in (\ref{Toeplitz}) is the inner product in ${L}^{2}(\CC, d\mu)$ defining the Bargmann projection (\ref{projector}), then
\beq
(T_\varphi\psi)(z) = (P ~ \varphi \psi) (z).   %( K_{\overline{z}} ,\varphi \psi ).
\label{ToeplitzBargmann}
\eeq
This last equation {\em defines} a Toeplitz operator on $F^2(\CC, d\mu)$ (see for instance \cite{Hal}).

According to original works by Berezin \cite{B1}, for a given symbol
$\varphi$ the \emph{natural domain} of  $T_\varphi$ is the subset of
${F}^{2}(\CC, d\mu)$ 
\beq 
\textit{Dom}(T_{\varphi}) = \{ \psi \in
{F}^{2}(\CC, d\mu): \, \varphi \psi \in L^{2}(\CC, d\mu) \}.
\label{B_domain} 
\eeq 
Then the Bargmann projection is well defined,
warranting that the integral in eq.\ (\ref{Toeplitz}) is an inner product 
in ${L}^{2}(\CC, d\mu)$, with result in  ${F}^{2}(\CC, d\mu)$.
Following this line, one can analyze classes of symbols such that
the natural domain is wide enough in  ${F}^{2}(\CC, d\mu)$. For
instance Berger and Coburn \cite{BergerCoburn} developed a symbolic
calculus for Toeplitz operators with bounded symbols, whose natural domains are the whole space ${F}^{2}(\CC, d\mu)$. 
Later, Coburn \cite{Coburn} defined a class of symbols $\varphi$ such that
\beq
\forall z \in \CC ~~~~ \varphi(w,\overline{w})e^{w\overline{z}} \, \in \, L^{2}(\CC, d\mu)\, ;
\label{Coburn_class}
\eeq
in other words, the Segal-Bargman symbols for all the coherent states belong to the natural
domain, making it dense in ${F}^{2}(\CC, d\mu)$.
In the same spirit, Folland \cite{Folland} considers safe to work with a class of symbols $\varphi$ such that
\beq
|\varphi(w,\overline{w})| \leq C \exp(\delta|w|^2), \text{with $\delta<1/2$, }
\label{Folland_class}
\eeq
which are included in Coburn's class.

However, the standard domain definition (\ref{B_domain}) is somewhat restrictive, as eq.\ (\ref{Toeplitz})
may make sense as a well defined integral even when it is not an inner product in $L^{2}(\CC;d\mu)$.
Some authors, as Janas \cite{Janas1991}, consider the largest domain in ${F}^{2}(\CC, d\mu)$:
\beq
\{\psi \in {F}^{2}(\CC, d\mu): \,  \int_\CC d\mu(w)\overline{ K(w,\overline{z})} \varphi(w,\overline{w}) \psi(w)
\in F^{2}(\CC, d\mu )\},
\label{J_domain}
\eeq
which in some cases may be indeed larger than $\textit{Dom}(T_{\varphi})$ \cite{Janas1994}.

% In this framework, let us finally recall that precise properties of Toeplitz operators discussed in the literature 
% heavily rely on the class of symbols involved. 
In what follows, we adhere to the natural domain definition (\ref{B_domain}) given by Berezin.

Besides the realization in the Segal-Bargmann space $F^2(\CC, d\mu)$,
for a given symbol $\varphi(w,\overline{w})$ one could also consider the realization of
the anti-Wick operator (\ref{antiWick}) on the other isometric spaces.
It is well known that the realization in $L^2(\RR, dx)$ leads to pseudo-differential operators in Weyl form \cite{Guillemin},
a setting where powerful tools are available.
The realization on $l^2$ has also been considered, in particular in \cite{G-V}, 
and turns out to be most convenient for Toeplitz operators with radial symbols.
It presents some features that will be analyzed in Sections \ref{V} and \ref{VI}.

%%%%%%%%%%%%%%%%%%%%%%%%%%%%%%%%%%%%%%%%%%%%%%%%%%%%%%%%%%%%%%%%%%%%%%%%%%%%%%%%%%%%%%%%%%%%%%%%%%%%%%%%%%%%%%%%%%%%%%%%%%

\section{Radial anti-Wick and Toeplitz operators} \label{IV}

We consider in this Section anti-Wick or Toeplitz operators with radial symbols, namely
$\varphi(w,\overline{w}) =\varphi(|w|)$. 
A salient feature of such operators is that any vector in the canonical basis (\ref{Fockbasis}), 
also belonging to the operator domain, is an eigenvector.
In order to discuss this property and its consequences, 
we find it convenient to give the following 

\begin{definition}
We denote ${\mathcal P}$ the class of radial symbols
\beq
{\mathcal P} = \left\{ \varphi(|w|): \forall n \in \NN, \, u_n \in \textit{Dom}(T_{\varphi}) \right\}
\label{Pclass}
\eeq
where $u_n(z)=  z^{n}/ \sqrt{n!}$ are the Segal-Bargmann symbols for vectors
$|n\ket$ in the orthonormal complete set (\ref{Fockbasis}), 
and $\textit{Dom}(T_{\varphi})$ is the natural domain defined in (\ref{B_domain}).
\label{Pdefined}
\end{definition}

\vspace{5mm}
For symbols $\varphi(|w|)$ in this class, the Toeplitz operator $T_{\varphi}$ has a dense natural domain, 
including at least any polynomial in $F^2(\CC, d\mu)$,
as well as $A_{\varphi}$ is well defined on any finite linear combination of vectors $|n\ket$.

In connection with the classes (\ref{Coburn_class}, \ref{Folland_class}) described in Section \ref{III}, we mention the following 
\begin{pro}
\label{PinCoburn}
If  $\varphi(|w|)$  is a radial symbol in the Coburn's class (\ref{Coburn_class}),
then $\varphi(|w|) \in {\mathcal P}$
\end{pro}

{\em Proof:} This relation is proven in Theorem 1.3 (i) in \cite{Janas1994}, even for non-radial symbols.
It is also proven that the converse is not true.

\vspace{5mm}
Let $\varphi(|w|) \in {\mathcal P}$ and consider the radial anti-Wick operator 
\beq 
A_\varphi = \int_\CC \frac{d^2w}{\pi}
\, |\overline{w}\ket \varphi(|w|) \bra \overline{w} | 
\label{radialAntiWick} 
\eeq 
Then the matrix elements $\bra m |A_\varphi| n \ket$ in the canonical basis (\ref{Fockbasis}) can be computed.
The key feature of radial symbols is that the integrals over the complex plane are easily solved in polar coordinates:
using eq.\ (\ref{zn}) it is straightforward to show that 
\beq 
\bra m |A_\varphi| n \ket 
= 
\int_\CC d\mu(w) \varphi(|w|)\frac{\overline{w}^m}{\sqrt{m!}}\frac{w^n}{\sqrt{n!}} 
=
\delta_{mn} \varphi_n,
\label{radiagonal} 
\eeq 
where % $\delta_{mn}$ is the Kroenecker delta and 
\beq
\varphi_n = \frac{2}{n!} \int_0^\infty \varphi(r) r^{2n+1} e^{-r^2} dr.
\label{phin0}
\eeq 
As the off-diagonal elements vanish the anti-Wick operator is drastically simplified, with a diagonal expression  
\beq 
A_\varphi =\sum_{n \in \NN} | n \ket \varphi_n \bra n |
\label{A_diag} 
\eeq
in the canonical Fock space basis (\ref{Fockbasis}). 
This of course provides its spectral decomposition, with eigenvalues $\varphi_n$ and eigenvectors $|n\ket$ in (\ref{Fockbasis}).

From this decomposition, it is most suited to consider the realization of $A_\varphi$ on $l^2$:
given a Fock space vector $| \psi \ket = \sum_{n\in\NN} \psi_n |n\ket$, one gets
\beq
A_\varphi| \psi \ket = \sum_{n\in\NN}\varphi_n \psi_n |n\ket.
\label{A_diag_psi}
\eeq
While  $| \psi \ket$ is realized by the sequence $\{\psi_n\}\in l^2$, its image $T_\varphi| \psi \ket$ 
is realized by the sequence  $\{\varphi_n \psi_n\}$. 
Introducing the notation $D_\varphi$ for the realization of $A_\varphi$ on $l^2$, we then have 
\beq
D_\varphi \{\psi_n\} = \{\varphi_n \psi_n\},
\label{D_l2}
\eeq
that is $D_\varphi$ acts as a pointwise multiplication operator on $l^2$ sequences.

Consider now the realization of the same anti-Wick operator, with symbol $\varphi(|w|) \in {\mathcal P}$ 
as a radial Toeplitz operator $T_\varphi$ on $F^2(\CC, d\mu)$.  
As a consequence of the isomorphism $U$ in Proposition \ref{prop_U}, 
it is apparent that the Toeplitz operator $T_\varphi$ on $F^2(\CC, d\mu)$ is unitarily equivalent 
to the diagonal operator $D_\varphi$ on $l^2$. 

\vspace{5mm}
The $l^2$ realization of radial Toeplitz operators has been used 
in \cite{Coburn, Coburn2010} and elsewhere to analyze simple examples, 
while it has been exploited by Grudsky and Vasilevski in \cite{G-V} as the key point 
to analyze boundedness, compactness and spectral properties of radial Toeplitz operators. 
However, the equivalence between Toeplitz operators with
radial symbols and diagonal operators of the form (\ref{D_l2}) must be treated carefully;
in particular, an unattentive reading of \cite{G-V} can be misleading, as operator domains are not explicitly 
regarded. 
To be precise, let us denote by $L_1^\infty(\RR_+,e^{-r^2})$, 
as in  \cite{G-V}, the set of all measurable functions $\phi(r)$ on $\RR_+$ such that 
\beq
\int_0^\infty |\phi(r)| r^{m} e^{-r^2} dr < \infty \, .
\label{L1}
\eeq
There exist symbols $\phi(|w|) \in L_1^\infty(\RR_+,e^{-r^2})$  but $\phi(|w|) \notin {\mathcal P}$ (see Example \ref{exists} below).
Though the matrix elements in (\ref{radiagonal}) can not be computed for such symbols,  
one still can compute a sequence $\{\phi_n\}$ as in (\ref{phin0}). We then give the following
\begin{definition}
\label{Dtilde_def}
Given a radial symbol $\phi(|w|)$ in $L_1^\infty(\RR_+,e^{-r^2})$, 
let $\tilde{D}_\phi: D \subset l^2 \to  l^2$ be defined by 
\beq
\tilde{D}_\phi \{\psi_n\} := \{\phi_n \psi_n\},
\label{Dtilde_l2}
\eeq
with 
\beq 
\phi_n := \frac{2}{n!} \int_0^\infty \phi(r) r^{2n+1} e^{-r^2} dr.
\label{phin}
\eeq 
\end{definition}
Notice that we define the domain $D$ of $\tilde{D}_\phi$ with the condition $\{\phi_n \psi_n\} \in l^2$. 
$D$ is dense in $l^2$, as $\tilde{D}_\phi$ is well defined at least on finite sequences, 
but for some infinite sequence $\{\phi_n \psi_n\}$ could lie outside $l^2$.

Now, if a radial symbol  $\phi(|w|) \in L_1^\infty(\RR_+,e^{-r^2})$  but $\phi(|w|) \notin {\mathcal P}$,
the operator $\tilde{D}_\phi$ has a dense domain in $l^2$
and then can not be equivalent to $T_\phi$.
This remark should be contrasted with Theorem 3.1 in \cite{G-V}. 
In order to complete the analysis given there, we summarize the above discussion in the following:

\begin{theorem}
Given a radial symbol $\varphi(|w|) \in L_1^\infty(\RR_+,e^{-r^2})$, 
the Toeplitz operator $T_\varphi$ defined in (\ref{ToeplitzBargmann})
is unitarily equivalent to the diagonal operator $\tilde{D}_\varphi$ defined in (\ref{Dtilde_def}) if and only if
$\varphi(|w|) \in {\mathcal P}$. In the positive case, the equivalence is given by $T_\varphi = U^{-1} \tilde{D}_\varphi U$, 
with the unitary operator $U$ defined in (\ref{U}).
\label{main}
\end{theorem}

{\em Proof}: 
If $\varphi(|w|) \in {\mathcal P}$, 
it is enough to consider the action of $U^{-1} \tilde{D}_\varphi U$ on functions $u_m$ as given in the definition \ref{Pdefined}.
$$
(U^{-1} \tilde{D}_\varphi U \, u_m)(z) = (U^{-1} \tilde{D}_\varphi)(\{\delta_{mn}\}_{n\in\NN})
=U^{-1}(\{\varphi_m \delta_{mn}\}_{n\in\NN})=\varphi_m u_m(z)
$$
On the other hand, it is straightforward to compute $(T_\varphi u_m)(z)$ using polar coordinates, getting the same result
$$
(T_\varphi u_m)(z)=\varphi_m u_m(z).
$$
Turning to the case where  $\varphi(|w|) \notin {\mathcal P}$, there exists $u_m \notin \textit{Dom}(T_{\varphi})$. 
But $U^{-1} \tilde{D}_\varphi U$ is well defined on $u_m$, then $T_\varphi$ and $U^{-1} \tilde{D}_\varphi U$ 
have different domains and in consequence are not equal.
\begin{flushright}
$\square$
\end{flushright}

It is important to show that the negative situation in Theorem \ref{main} is indeed possible with the following:

\begin{example}  
\label{exists}
Consider the function $\varphi(r)=e^{(\frac{1}{2}+\frac{\sqrt{3}}{2}i)r^{2}}$.
It belongs to \\
$L_1^\infty(\RR_+,e^{-r^2})$, with an associated sequence 
$\varphi_n = (\frac{1}{2}-\frac{\sqrt{3}}{2}i)^{-(n+1)}$. 
Moreover, this sequence is bounded ($|\varphi_n|=1$) and the operator (\ref{Dtilde_l2}) is bounded, with domain $l^2$.
However, $\varphi(|w|) \notin {\mathcal P}$: it is enough to observe that $u_0(z)=1$ does not belong to $\textit{Dom}(T_{\varphi})$.
\end{example}

\vspace{5mm}
We stress that, if $\varphi(|w|) \in {\mathcal P}$, then $\tilde{D}_{\varphi}$
is nothing but the realization on $l^2$ of the formal anti-Wick operator $A_{\varphi}$
realized by the Toeplitz operator $T_{\varphi}$ on $F^2(\CC, d\mu)$. 
On the other hand, if $\varphi(|w|) \notin {\mathcal P}$, one can consider the following 

\begin{definition}
\label{Ttilde-def}
Given a radial symbol $\varphi(|w|) \in L_1^\infty(\RR_+,e^{-r^2})$, 
let \\
$\tilde{T}_\varphi: D \subseteq F^2(\CC, d\mu) \to  F^2(\CC, d\mu)$ be defined by 
\beq
\tilde{T}_{\varphi } := U^{-1} \tilde{D}_\varphi U
\label{Ttilde-eq}
\eeq
with $\tilde{D}_\varphi$ given in \ref{Dtilde_def} and  $U$ defined in (\ref{U}),
with domain \\
$D= \textit{Dom}(\tilde{T}_{\varphi})=\left\{ f \in F^2(\CC, d\mu): \tilde{D}_\varphi U f \in  l^2 \right\}$.
\end{definition}

In order to motivate this definition we prove the following:

\begin{theorem}
Let $T_\varphi$ be a radial Toeplitz operator with symbol $\varphi$ in\\
 $L_1^\infty(\RR_+,e^{-r^2})$,
and let $\tilde{T}_{\varphi}$ be defined in \ref{Ttilde-def}. Then\\ 
a - $\textit{Dom}(T_{\varphi}) \subseteq \textit{Dom}(\tilde{T}_{\varphi})$,  and \\
b - the restriction of $\tilde{T}_{\varphi}$ to $\textit{Dom}(T_{\varphi}) $ coincides with $T_\varphi$.
\label{extension}
\end{theorem}

{\em Proof}:

Let $\varphi \in L_1^\infty(\RR_+,e^{-r^2})$, and let $f(z) \in \textit{Dom}(T_{\varphi})$.
As $f(z) \in F^2(\CC, d\mu)$, it can be expanded as $f(z) = \sum_{n \in \NN} a_n u_n(z)$ with $\{a_n\} \in l^2$.
Let $g(z) = (T_{\varphi} f)(z)$, also in $F^2(\CC, d\mu)$, be expanded as $g(z) = \sum_{n \in \NN} b_n u_n(z)$ with $\{b_n\} \in l^2$.
The coefficients $b_n$ can be computed with the inner product in $F^2(\CC, d\mu)$ as $b_n=(u_n,g)$; 
using polar coordinates one easily obtains
\beq
b_n=(u_n,g)=\varphi_n a_n,
\label{gn}
\eeq
showing that $\{\varphi_n a_n\} \in l^2$.

a - To see that $f \in \textit{Dom}(\tilde{T})$, compute
$$
\tilde{T}_{\varphi} f = U^{-1} \tilde{D}_\varphi U f =  U^{-1} \tilde{D}_\varphi \{a_n\} = U^{-1} \{\varphi_n a_n\}.
$$
As $\{\varphi_n a_n\} \in l^2$, then $f \in \textit{Dom}(\tilde{T}_{\varphi})$. Moreover,
$$
\tilde{T}_{\varphi} f  = \sum_{n \in \NN}\varphi_n a_n u_n(z) \in F^2(\CC, d\mu).
$$ 

b - For $f(z) \in \textit{Dom}(T_{\varphi})$ we have computed above $ (T_{\varphi} f)(z)$, coinciding with the last equation.

\begin{flushright}
$\square$
\end{flushright}

From the discussion above, given a radial symbol $\varphi(|w|) \in L_1^\infty(\RR_+,e^{-r^2})$
but $\varphi(|w|) \notin {\mathcal P}$, one can propose $\tilde{T}_\varphi$ to be considered a nontrivial extension 
of the Toeplitz operator $T_{\varphi}$, as in this case the inclusion 
$\textit{Dom}(T_{\varphi}) \subset \textit{Dom}(\tilde{T}_{\varphi})$ is strict.
To close this Section, we would like to point out that most of the results in \cite{G-V} apply to such an
extension.

%%%%%%%%%%%%%%%%%%%%%%%%%%%%%%%%%%%%%%%%%%%%%%%%%%%%%%%%%%%%%%%%%%%%%%%%%%%%%%%%%%%%%%%%%%%%%%%%%%%%%%%%%%%%%%%%%%%%%%%%%%

\section{Diagonal operators on $l^2$} \label{V}

As a radial Toeplitz operator can, under certain assumptions, be faithfully mapped into a diagonal operator on $l^2$, 
the natural question to consider is the inverse situation. 
Namely, given a sequence $\{\gamma_n\}$ defining the spectrum of a diagonal operator on the canonical basis of $l^2$, 
we would like to know whether there exists a radial symbol $\gamma(|w|)$ defining an equivalent Toeplitz operator.
% A similar question has been answered for Toeplitz operators acting on the Bergman space \cite{Suarez}.
But in the Segal-Bargmann space the answer to this question should as a first step invert eq.\ (\ref{phin}) to construct the symbol $\gamma(|w|)$, 
and, as a second step, establish whether this symbol provides a Toeplitz operator defined on any polynomial or not, 
as discussed in the previous Section.

The first step was partially solved in \cite{G-V}: by means of a nice construction based on analytic continuation and Fourier transforms, 
the authors show that for any {\em bounded} sequence $\{\gamma_{n}\}$ there exists a radial symbol $\gamma(|w|)$ in 
$L_1^\infty(\RR_+,e^{-r^2})$ such that the sequence defined in (\ref{phin}) is again $\{\gamma_n\}$. However, 
a symbol in this set does not guarantee the equivalence of a Toeplitz operator $T_\gamma$ with a diagonal operator $\tilde{D}_\gamma$.
In an attempt to complete the second step we found that, unfortunately, the symbol provided in \cite{G-V} is hard to analyze in full generality.
Moreover, as far as we know there are no results for {\em unbounded} sequences.

In this section we aim to provide some insight on these questions, by restricting to a family of 
sequences $\{k^{-n}\}_{n\in\NN}$ for which a symbol $\gamma_k$ in $L_1^\infty(\RR_+,e^{-r^2})$ can be recognized with a close expression:

\begin{pro}
\label{V1}
Let $k \in \CC$. If $Re(k)>0$, the sequence $\{k^{-n}\}$ is obtained by (\ref{phin})) from the radial symbol $\gamma_k(|w|) =k e^{(1-k)|w|^2} \in L_1^\infty(\RR_+,e^{-r^2})$.
\end{pro}

{\em Proof}: The computation is straightforward.

\vspace{5mm}
It is simple to explore some properties of sequences in \ref{V1}:

\begin{pro}
The symbol $\gamma_k(|w|)$ is in the class ${\mathcal P}$ if and only if $Re(k)>1/2$. Correspondingly, 
the Toeplitz operator $T_{\gamma_k}$ is unitarily equivalent to $\tilde{D}_{\gamma_k}$ defined in \ref{Dtilde_def} if and only if $Re(k)>1/2$.
\label{V2}
\end{pro}

{\em Proof}: according to the definition in (\ref{Pclass}), it is easy to check that for all  $n \in \NN$
$$ \int d\mu(w) |\gamma_k(|w|) u_n(w)|^2 < \infty  $$
if and only if $Re(k)>1/2$. 
From Theorem \ref{main}, $T_{\gamma_k}$  is then unitarily equivalent to $\tilde{D}_{\gamma_k}$.

\vspace{5mm}
Relatedly, the same condition on $k$ characterizes whether the symbol $\gamma_k$ 
is in the classes considered by Folland and Coburn:

\begin{pro} ~
\\
a - The symbol $\gamma_k$ is in Folland's class (\ref{Folland_class})  if and only if $Re(k)>1/2$.\\
b - The symbol $\gamma_k$ is in Coburns's class (\ref{Coburn_class})  if and only if $Re(k)>1/2$.
\label{V3}
\end{pro}

{\em Proof}: \\
a - notice that $\gamma_k(|w|)=k  e^{i Im(k) |w|^2} e^{\delta |w|^2} $
with $\delta = 1 - Re(k)$.\\
b - as in the previous proposition, it is easy to check that for all  $z \in \CC$
$$ \int d\mu(w) |\gamma_k(|w|) K_{z}(w)|^2 < \infty  $$
if and only if $Re(k)>1/2$. 

\vspace{5mm}
As a consequence, our symbols $\gamma_k(|w|)$ belong simultaneously to the classes ${\mathcal P}$,
Folland's and Coburn's, or to none of them.  While no subtle differences between this classes will be found in our analysis,
the results in this Section apply to all of them.
Notice that a well known operator in Quantum Mechanics, the Maxwell-Boltzman density matrix for a harmonic oscillator, 
can be written as a Toeplitz operator with this kind of symbol.
The operator is defined as
$$\widehat\rho = e^{-\beta ( \widehat{a}^{*} \widehat{a} +1/2)},$$ 
where $\beta > 0$ is the (dimensionless) inverse temperature.
As 
$$\widehat\rho |n\ket> = e^{-\beta/2}e^{-\beta n}|n\ket,$$
the operator is diagonal on $l^2$ with a power like sequence spectrum \\
$\{e^{-\beta/2} (e^\beta)^{-n}\}$.
Being $e^\beta > 1$, the corresponding anti-Wick symbol reads $\rho(|w|) = e^{\beta/2}e^{(1-e^\beta)|w|^2}$
and belongs to ${\mathcal P}$.

We explore first bounded sequences $\{k^{-n}\}$, with $|k| \geq 1 $.

\begin{pro}
If $Re(k)>1/2$ and $|k| \geq 1$, the symbol $\gamma_k(|w|)$ coincides with the one constructed by 
Grudsky and Vasilevski in the Theorem 3.7 of \cite{G-V}.
Moreover, the Toeplitz operator $T_{\gamma_k}$ is then unitarily equivalent to $\tilde{D}_{\gamma_k}$.
\label{V4}
\end{pro}

{\em Proof}: as $|k| \geq 1$, the sequence $\{k^{-n}\}$ is bounded and Theorem 3.7 of \cite{G-V} 
provides a radial symbol in $L_1^\infty(\RR_+,e^{-r^2})$. As $Re(k)>1/2$,
uniqueness of the symbol for a given Toeplitz operator follows as a known property of Folland's class \cite{Folland}.
Also equivalence between $T_{\gamma_k}$ and  $\tilde{D}_{\gamma_k}$ follows from Theorem \ref{main}.

\begin{pro}
If $|k|>1$ and $0 < Re(k) < 1/2$, the symbol $\gamma_k$ belongs to $ L_1^\infty(\RR_+,e^{-r^2})$
but the natural domain $\textit{Dom}(T_{\gamma_k})$ is trivial. 
\label{V5}
\end{pro}

{\em Proof}: As mentioned in Proposition \ref{V1}, $ Re(k) > 0$ implies that $\gamma_k(|w|) \in  L_1^\infty(\RR_+,e^{-r^2})$.
For $f \in F^2(\CC, d\mu)$ being in the natural domain of $T_{\gamma_k}$, the condition $\gamma_k f \in L^2(\CC, d\mu)$
requires $|f(z)|$ to decay exponentially; then, by Liouville theorem, the only function in  $\textit{Dom}(T_{\gamma_k})$ is
$f(z)=0$.

\vspace{5mm}
Turning to unbounded sequences, $|k|<1$, there is no general statement about the existence of an equivalent Toeplitz operator.
However, for the family of symbols under consideration, we learn that:

\begin{pro}
If   $|k|<1$ and $ Re(k) > 1/2$, the sequence $\{k^{-n}\}$ is unbounded  but the symbol $\gamma_k \in {\mathcal P}$ 
still allows the definition of an associated $\tilde{D}_{\gamma_k}$ as defined in \ref{Dtilde_def}.
Then the Toeplitz operator $T_{\gamma_k}$ is densely defined, unbounded and unitarily equivalent to 
the diagonal operator $\tilde{D}_{\gamma_k}$ on $l^2$.
\label{V6}
\end{pro}

% The following table provides explicit examples of the different situations discussed above:
% 
% \begin{center}
% \begin{tabular}{|c|c|c|c|}
%   \hline
% symbol $\gamma$  & sequence & boundedness & $\varphi \in {\mathcal P}$ \\
%   \hline
%   % after \\: \hline or \cline{col1-col2} \cline{col3-col4} ...
%   $e^{1/2r^{2}}$ & $2^{n+1}$ & no & no \\
%   $e^{1/3r^{2}}$ & $(3/2)^{n+1}$ & no & yes \\
%   $e^{(1/2+\sqrt{3}/2)r^{2}}$ & $(1/2+\sqrt{3}/2)^{-(n+1)}$ & yes & no \\
%  $e^{(1/3+\sqrt{8}/3)r^{2}}$ & $(1/3+\sqrt{8}/3)^{-(n+1)}$ & yes & yes \\
%   \hline
% \end{tabular}
% \end{center}

%%%%%%%%%%%%%%%%%%%%%%%%%%%%%%%%%%%%%%%%%%%%%%%%%%%%%%%%%%%%%%%%%%%%%%%%%%%%%%%%%%%%%%%%%%%%%%%%%%%%%%%%%%%%%%%%%%%%%%%%%%

\section{Composition of Toeplitz
operators: an open problem.} \label{VI}

A problem still open is to determine when the composition of Toeplitz
operators is closed. That is, given symbols $\varphi$ and $\eta$ in
certain class, to determine the existence of a symbol $\tau$ such
$T_{\tau}$ is a well defined Toeplitz operator and $T_{\tau} =
T_{\varphi} T_{\eta}$.

A first approach to this question was given by Coburn in \cite{Co},
where symbols in a smooth Bochner algebra are considered. These
are Fourier-Stieltjes transforms of compactly supported, regular,
bounded complex-valued Borel measures on $\CC$, indeed a subclass of
the symbols defined in (\ref{Coburn_class}). The author proved that this class 
is closed under composition: $T_{\varphi} T_{\eta}$ is a Toeplitz
operator with a symbol $\tau$ also in the smooth Bochner algebra.
The composition symbol $\tau$ can be calculated as a Moyal type
product $\tau= \varphi \diamond \eta$ given by
\begin{equation}
(\varphi  \diamond \eta)(w,\overline{w}) =
\sum_{k}\frac{(-1)^{k}}{k!}\partial^{k}\varphi(w,\overline{w})\overline{\partial}^{k}\eta(w,\overline{w}).
\label{diamond}
\end{equation}
where $\partial^{k}:=\frac{\partial^{k}}{\partial^{k}w}$ and
$\overline{\partial}^{k}:=\frac{\partial^{k}}{\partial^{k}\overline{w}}$.
Later, in \cite{Coburn} this result was extended to unbounded
Toeplitz operators with polynomial symbols in $w$ and
$\overline{w}$.

More recently Bauer \cite{Bauer} considered in detail the class of
symbols in the set 
\beq Sym_{>0}(\CC):=
\bigcap_{j=1}^{\infty}\mathcal{D}_{1/j} \label{Bauer_class} 
\eeq
where 
$\mathcal{D}_c := \{ \varphi : \exists d >0 \
\mbox{such that} |\varphi(z,\overline{z})| \leq \ d \ exp(c|z|^2) \
\ \ a.e. \}$, 
providing well defined unbounded Toeplitz operators
with domain on a certain scale of Banach spaces. The author proved
that the composition of Toeplitz operators is closed for this class, 
with a symbol for the composition given by the same product (\ref{diamond}).
Notice that polynomial symbols are included here providing an
independent proof of Coburn's results.

For Toeplitz operators with radial symbols, the composition problem
can be simplified using the results in Section \ref{IV}. 
Given two radial symbols $\varphi$, $\eta$ in the class ${\mathcal P}$,
from Theorem \ref{main} one has 
\beq
T_{\varphi}= U^{-1} \tilde{D}_{\varphi} U
\ \ \ \ \ \ \ \ \ \ \ \ \ \ \ \
T_{\eta}= U^{-1} \tilde{D}_{\eta} U
\label{toep->diag} 
\eeq
so that 
\beq
T_{\varphi} T_{\eta}= U^{-1} \tilde{D}_{\varphi} \tilde{D}_{\eta} U,
\label{compose-toep} 
\eeq
where the composition of diagonal operators on $l^2$ is trivially expressed as 
\beq
\tilde{D}_{\varphi} \tilde{D}_{\eta}\{c_n\} = \{\varphi_n \eta_n c_n\}.
\label{compose-diag} 
\eeq
The domain of this operator is clearly dense in $l^2$, so that 
the composition $ T_{\varphi} T_{\eta} $ in (\ref{compose-toep}) is well defined on a dense domain in $F^2(\CC, d\mu)$.
It is also clear that the composition is commutative in this case.
However, $ U^{-1} \tilde{D}_{\varphi} \tilde{D}_{\eta} U $ does not necessarily correspond to a radial Toeplitz
operator with some symbol $\tau$. From the results in Section \ref{IV}, can summarize the following conditions:

\begin{pro}
\label{P-closed}
Let $\varphi$, $\eta$ be radial symbols in the class ${\mathcal P}$.
If there exists a radial symbol $\tau$ in ${\mathcal P}$ such that the sequence $\{\tau_n\}$ 
given by (\ref{phin}) coincides with $\{\varphi_n \eta_n \}$, 
then $T_{\tau} = T_{\varphi} T_{\eta}$, that is the composition of radial Toeplitz 
operators is another radial Toeplitz operator. 
If not, the composition of such radial Toeplitz operators is not a Toeplitz operator with symbol in ${\mathcal P}$.
\end{pro}

There are known counter examples \cite{Bauer}-\cite{Coburn} that
exhibit such limitations on the possibility of composing Toeplitz
operators with symbols in some given class. 
In particular, \cite{Coburn} discusses an example with a radial symbol which is indeed 
in the family that we considered in Section \ref{V} (with $k=\frac{3}{5}-\frac{4}{5}i$).
We can give some generality to this analysis, within the symbols in that family, with the following: 
\begin{pro}
\label{compose}
Let $\gamma_a(|w|)=a e^{(1-a)|w|^{2}}$ and $\gamma_b(|w|)=b e^{(1-b)|w|^{2}}$ with
$Re(a)>1/2 $ and $Re(b)>1/2$. The operator $T_{\gamma_a} T_{\gamma_b}$ is
a Toeplitz operator with radial symbol in ${\mathcal P}$ if and only if $Re(ab)>1/2$. 
In the positive case, the symbol is $\gamma_c= c e^{(1-c)|w|^{2}}$ with $c=ab$.
\end{pro}

{\em Proof:} From \ref{V1} and  \ref{V2}, $T_{\varphi}$ and $T_{\eta}$ are unitarily equivalent to 
$\tilde{D}_{\varphi}$ and $\tilde{D}_{\eta}$, described by sequences $\{a^{-n}\}$ and $\{b^{-n}\}$, respectively. 
Their composition is described by the sequence $\{(ab)^{-n}\}$. 
The condition on $Re(ab)$ and the symbol $\gamma_c$ are finally read from \ref{V1} and \ref{V2}.

\vspace{5mm}
It is interesting to note that: 

\begin{pro}
\label{Moyal}
The Moyal product (\ref{diamond}) given by Coburn for polynomial and Bochner algebra symbols   
provides the correct symbol for the composition described in Proposition \ref{compose}, with $Re(a.b)>1/2$.
\end{pro}

{\em Proof:} The computation of the Moyal product is straightforward, and the resulting series can be summed 
to give the expected result. 
Indeed, the product can be done for any $a,b \in \CC$ and is commutative, but when the symbols are not in ${\mathcal P}$ 
one can not associate densely defined Toeplitz operators to them.  

\vspace{5mm}
The composition of Toeplitz operators is then not closed for the family of symbols $k e^{(1-k)|w|^{2}} \in {\mathcal P}$ with $Re(k) > 1/2$. 
However, particular examples where composition of 
Toeplitz operators not belonging to closed classes is a Toeplitz operator can be given.
For instance, pick $a=\frac{3}{5}-\frac{4}{5}i$ and $b=\frac{3}{5}+\frac{4}{5}i$, so that $ab=1$,
$T_{\gamma_a} T_{\gamma_b} =T_{\gamma_1 }$ is the identity in $F^2(\CC, d\mu)$
and all the involved symbols are in ${\mathcal P}$.

One is tempted now to look for a class of radial symbols where the composition of the corresponding Toeplitz operators is closed.
Within the scope of Proposition \ref{V1}, we are forced to restrict to the following:

\begin{definition}
Let 
$$
{\mathcal L} = \left\{ \gamma_k(|w|)=k e^{(1-k)|w|^{2}}:\, k\in\RR \wedge k \geq 1 \right\} .
$$
\label{small_class}
\end{definition}

From Proposition \ref{compose}, it is immediate that:
\begin{pro}
The class of radial symbols $ {\mathcal L}$ is closed under composition of the corresponding Toeplitz operators.
\end{pro}
Unfortunately, this result is not new, as symbols in $ {\mathcal L}$ are bounded. 
They are in Bauer's class defined in (\ref{Bauer_class}), and also
in the (not smooth) Bochner algebra, being the Fourier-Stieltjes of Gaussian measures with 
non compact support.

We conclude this section reviewing the mentioned example in \cite{Coburn},
the Toeplitz operator $T_{\gamma_a}$ with  $a=\frac{3}{5}-\frac{4}{5}i$. 
The symbol $ \gamma_a(|w|)$ is not bounded, but
the natural domain of  $T_{\gamma_a}$ contains all of the coherent states, 
that is $ \gamma_a(|w|)$ belongs to the class defined in (\ref{Coburn_class}).
Coburn proves that the composition $T_{\gamma_a} T_{\gamma_a}$ can not be written as a Toeplitz operator with a symbol in the same class.
In our notation, as $Re(a)>1/2$, we can write
$$ 
T_{\gamma_a} = U \tilde{D}_{\gamma_a} U^{-1},
$$
where
$
\tilde{D}_{\gamma_a}\{\psi_n\} = \{a^{-n} \psi_n \}.
$
As $|a|=1$, the spectrum of $\tilde{D}_{\gamma_a}$ (equivalently $T_{\gamma_a}$) is bounded. 
Then the natural domain of $T_{\gamma_a}$ is the full space $F^2(\CC, d\mu)$. 
The composition of $T_{\gamma_a}$ with itself reads
$$ 
T_{\gamma_a} T_{\gamma_a} = U \tilde{D}_{\gamma_a} \tilde{D}_{\gamma_a} U^{-1},
$$
well defined on $F^2(\CC, d\mu)$ as 
$
\tilde{D}_{\gamma_a} \tilde{D}_{\gamma_a} \{\psi_n\} = \{(a^2)^{-n} \psi_n \}
$
has again bounded spectrum.
However, $a^2=\frac{-7}{25}-\frac{24}{25}i$ has negative real part:
according to Proposition \ref{compose}, $T_{\gamma_a} T_{\gamma_a}$ is not a Toeplitz operator with symbol in ${\mathcal P}$,
a result that in our framework corresponds to that proved in \cite{Coburn}.

Notice that an operator $\G:l^2 \to l^2$ defined as $\G:=\tilde{D}_{\gamma_a} \tilde{D}_{\gamma_a}$ 
is diagonal, acting as pointwise multiplication by the bounded sequence $\{(a^2)^{-n}\}$. 
Then, according to Theorem 3.7 in \cite{G-V},
there exists a radial symbol $\varphi(|w|)$ in $L_1^\infty(\RR_+,e^{-r^2})$ such that eq.\ (\ref{phin}) 
generates the sequence $\{(a^2)^{-n}\}$. This symbol can not belong to ${\mathcal P}$, so that
a Toeplitz operator $T_\varphi$ with this symbol is not densely defined in $F^2(\CC, d\mu)$
and is not equivalent to $T_{\gamma_a} T_{\gamma_a}$.

%%%%%%%%%%%%%%%%%%%%%%%%%%%%%%%%%%%%%%%%%%%%%%%%%%%%%%%%%%%%%%%%%%%%%%%%%%%%%%%%%%%%%%%%%%%%%%%%%%%%%%%%%%%%%%%%%%%%%%%%%%

\vspace{1cm}
{\em Acknowledgments}: the authors wish to thank G. Silva for helpful discussions. 
This work was partially supported by CONICET (grant PIP 1691) and ANPCyT (grant PICT 20350), Argentina.

\vspace{2cm}

\section*{Appendix: the standard coherent states set}

The set of coherent states in eq.\ (\ref{CS})
is not an orthogonal basis but an overcomplete set in ${\mathcal F}$ \cite{Perelomov}.
Indeed, one can show that
\beq
\bra\overline{z} | w \ket = e^{-\frac{1}{2}|z|^2} e^{-\frac{1}{2}|w|^2} e^{z w} \, ,
\label{zw}
\eeq
while the completeness is given by a resolution identity in ${\mathcal F}$, {\em i.e.}
\beq
\mathbb{I} = \int_\CC \frac{\d2z}{\pi} | z \ket\bra z |
\label{identity}
\eeq
with the integral understood in the weak sense.
The notation $\bra z |$ in eq.\ (\ref{identity}) stands for the linear form $\bra z | \cdotp \ket$ on ${\mathcal F}$.

Vectors $|z\ket$ are conveniently related to the vacuum by
\beq
|z\ket = e^{-\frac{1}{2}|z|^2}e^{z\widehat{a}^{*}} |0\ket .
\label{z_def}
\eeq
Using eq.\ (\ref{z_def}), the operators $\widehat{a}^{*}$ and $\widehat{a}$ can be shown to be realized on  $F^{2}(\CC, d\mu)$
functions by
$$
a^{*} \psi(z) = z \psi(z), ~~~~~~~~ a \psi(z) = \frac{\partial}{\partial z} \psi(z).
$$
Then, vectors $|n\ket$ in the orthonormal basis (\ref{Fockbasis}) have Segal-Bargmann symbols
\beq
u_n(z) = z^{n}/ \sqrt{n!}.
\label{n_symbol}
\eeq
The symbol for a normalized coherent state $|w\ket$ reads, according to eq.\ (\ref{zw}),
\beq
u_w(z) =  e^{-\frac{1}{2}|w|^2} e^{w z} .
\label{z_symbol}
\eeq

It follows from eq.\ (\ref{identity}) that
\beq
\psi(z)=\int_\CC d\mu(w) \overline{  K_{\overline{z}}(w)} \psi(w),
\label{F2id}
\eeq
with $K_{\overline{z}}(w)= e^{\overline{z} w }$ as in eq.\ (\ref{kernel}).
The identity in $F^{2}(\CC, d\mu)$ is then represented as an integral operator with kernel $K_{\overline{z}}(w)$,
or as an inner product in $F^{2}(\CC, d\mu)$,
\beq
\psi(z)= ( K_{\overline{z}}, \psi )\, .
\label{RKHS}
\eeq
In other words, the restriction of the Bargmann projection (\ref{projector})
to $F^{2}(\CC, d\mu)$ is the identity operator.
Thus the existence of the set $ \{K_z \}_{z \in \CC} \subseteq F^{2}(\CC, d\mu)$ makes it a reproducing kernel Hilbert space.
Notice in pass that vectors $K_{\overline{z}}(w)$ are proportional to coherent state vectors
$u_{\overline{z}}(w)$ in eq.\ (\ref{z_symbol})
(some authors prefer to disregard normalization and refer to $K_{\overline{z}}(w)$ as Poisson vectors):
eq.\ (\ref{RKHS}) can be read as a weighted orthogonal projection of $| \psi \ket$ on $| \overline{z} \ket$.

All of the above can be trivially generalized to a finite number $N$ of commuting creation and annihilation
pairs, leading to $F^{2}(\CC^N, d\mu)$ spaces. The case of infinite number of coordinates was considered by Segal \cite{Segal}.

\end{document}